\documentstyle[leqno,amssymb,amstex,
11pt]{amsart} 
\def 	     \n 		{\bf N} 
\def	     \r 		{\bf R} 
\newcommand {\usumgxi} 		{\sum_{i,j=1}^d g^{ij}\xi_i\xi_j}

\newcommand {\usumaxi} 		{\sum_{i,j=1}^d a^{ij}\xi_i\xi_j}
\newcommand {\lsumgxi} 		{\sum_{i,j=1}^d g_{ij}\xi^i\xi^j}
\newcommand {\lsumghxi} 	{\sum_{i,j=1}^d g_{ij,h}\xi^i\xi^j}
\newcommand {\lsumaxi} 		{\sum_{i,j=1}^d a_{ij}\xi^i\xi^j}
\newcommand {\usumgdu} 		{\sum_{i,j=1}^d g^{ij}D_iuD_ju}
\newcommand {\usumghdu} 	{\sum_{i,j=1}^d g_h^{ij}D_iuD_ju}
\newcommand {\usumadu} 		{\sum_{i,j=1}^d a^{ij}D_iuD_ju}

\newcommand {\lsumadu}   	{\sum_{i,j=1}^d a^{ij}D_iuD_ju}
\newcommand {\huh}       	{{H^1(M_h,g_h)}}
\newcommand {\hu}		{{H^1(M,g)}}
\newcommand {\locvolh}		{\sqrt{\det g_h}}
\newcommand {\locvol}		{\sqrt{\det g}}
\newcommand {\riema}		{(a_{ij})_{i,j=1}^d}
\newcommand {\uriema}		{(a^{ij})_{i,j=1}^d}

\newcommand {\e}		{\varepsilon}
\newcommand {\ie}		{{\em i.e.,\/}}

\numberwithin{equation}{section}

\newtheorem {th}            	{Theorem}
\newtheorem {lemma}[equation]     	{Lemma}
\newtheorem {prop}[equation]	{Proposition}

\theoremstyle{definition}
\newtheorem{definition}[equation]{Definition}

\theoremstyle{remark}
\newtheorem{remark}[equation]		{Remark}
\newtheorem{ex}[equation]		{Example}
\newtheorem{corollary} [equation]	{Corollary}

\begin{document}

\title{Spectrum of compact manifolds with high genus } 
\author{Lino Notarantonio}
\address{Department of Theoretical Mathematics \\
The Weizmann Institute of Science\\
Rehovot 76100, Israel}
\email{notaran@@wisdom.weizmann.ac.il}  


\maketitle

\begin{abstract}
In this paper we study the behavior of the spectrum of a compact,
connected Riemannian manifold $(M,g)$ of dimension $d \ge 2$, when we
add an increasing number of 
increasingly small handles. No assumptions on any of the curvatures
are needed.
\end{abstract}

\maketitle 
\pagestyle{myheadings}
\pagenumbering{arabic}


\section*{Introduction}

\noindent
Let $( M , g ) $ be a compact  (connected, Riemannian) manifold (with
or without boundary) of dimension $d \ge 2$; let $g = (g_{ij})_{i,j=1}^d$
be its  metric. We shall  find convenient to work in the category of
{\em Lipschitz\/} manifolds, namely the coefficients of the metric
$g_{ij}$, $i,j=1\ldots, d$, are bounded, measurable  functions and
changes of coordinates are related by bi-lipschitzian diffeomorphisms.

We introduce  the spectrum of a Lipschitz manifold $(M,g)$ in
Definition~\ref{hone-spectrum}; if $(M,g)$ happens to be a 
smooth manifold, then the spectrum of $(M,g)$ is related with 
sequence of eigenvalues of the Laplace-Beltrami operator $-\Delta_g$
of $(M,g)$ (with ``natural'' boundary conditions) in an obvious way;
cf. Definition~\ref{nat-eigen} and Remark~\ref{remark-hone}.    

In this paper we shall be
interested in exhibiting topological perturbations  of $M$ which
affect the spectrum of $(M,g)$; more precisely we shall add an
increasing number of increasingly small  handles to $M$, call $(N_h,
g^{N_h})$ the {\em resulting manifold\/}, and study the spectrum of 
$(N_h , g^{N_h})$ as $h \to \infty$.  

Our main result is Theorem~\ref{limit-handles} where we prove that,
provided that the handles have a suitable girth, the spectrum of the
Laplace-Beltrami operator of $(N_h , g^{N_h})$ with ``natural''
boundary condition (cf. Definition~\ref{nat-eigen}) converges to the
spectrum of the operator induced by the quadratic form 
\begin{align}
{\cal D}_\infty(u) := 
{\cal D}_{M} ( u )  & + \frac{ \alpha }{ 2^d }\int_{ M } \bigl[
u(x) - u(T(x)) \bigr]^2 \mbox{Vol}_g(dx) + \notag\\ 
& + \lambda \int_{ M } u^2\mbox{Vol}_g(dx), 
\ \ \ u\in H^1(M, g); \tag{*}
\end{align}
\( T:M\rightarrow M\) is an involution which is induced by the
process of attaching handles to $M$. The convergence of the spectrum
mentioned above means convergence of the eigenvalues and of the
eigenspaces (Definition~\ref{conv-spaces})
generated by the corresponding eigenfunctions. It should
also be noted that the involution $T$ induces an orthogonal
decomposition of both \( L^2(M,g)\) and $H^1(M, g)$ into ``odd''  and
``even'' subspaces; therefore formula (*) above implies that the
process of attaching handles affects only the odd part of the
spectrum; cf. the discussion after Theorem~\ref{limit-handles}.

The proof of our main result relies on two general results, 
Theorem~\ref{var-comp} and Theorem~\ref{handles-comp}, 
which have a stronger analytical flavor and may have some interest in
their own. By using a general result by G.~Dal Maso, R.~Gulliver \&
U.~Mosco in \cite{dmgm} (cf. Proposition~\ref{representation} below),
we prove Theorem~\ref{limit-handles}  in the framework  of {\em
relaxed\/} Lipschitz manifolds (Definition~\ref{relax-mfld}), the
reason being that relaxed Lipschitz manifolds are better suited for
the analysis via $\Gamma$-convergence that will be used in this paper.

We notice that Theorem~\ref{limit-handles}, Theorem~\ref{var-comp} and
Theorem~\ref{handles-comp} generalize some results in \cite[\S\S
4,5]{dmgm}: There  those authors considered Lipschitz
manifolds-with-boundary topologically equivalent to bounded open sets
of ${\bf R}^d$, while our results deal with compact Lipschitz
manifolds, possibly without boundary.

Theorem~\ref{limit-handles} is proved in the case of $(M,g) = (S^d ,
g)$, the $d$-dimensional sphere of radius 1, with its standard metric
of constant sectional curvature. We have made
this choice so as to approach and present the problem in a hopefully
intuitive way. Even though the metric of the ``base'' manifold $M=S^d$
is smooth, we stress that the metric of the manifold which results from
attaching handles to $M$ is in general only Lipschitz; 
cf. Remark~\ref{handle-line-ex}  and also Remark~\ref{lip-mfld}.
We also stress that, adapting some arguments in \cite[\S\S
2,3]{dmgm}, it is possible to prove Theorem~\ref{limit-handles} in the
full generality of Lipschitz manifolds \( (M,g)\).

Theorem~\ref{var-comp} and Theorem~\ref{handles-comp} deal with
variational limits of sequences of perturbed Dirichlet
functionals. Theorem~\ref{var-comp} may be thought of as a compactness
result for Lipschitz  metrics, while Theorem~\ref{handles-comp} is a
compactness result for non-local perturbations of their associated
Dirichlet functionals.  

We point out that the convergence of  manifolds
treated in Theorem~\ref{var-comp} is quite different
from the convergence of manifolds as studied, among others, by
R.E.~Greene, M.~Gromov, S.~Peters and H.~Wu \cite{gromov},
\cite{greenewu}, \cite{peters} in that we need no bounds on any of the
curvatures; cf. Remark~\ref{handle-line-ex}. 

To obtain our result of convergence of the spectrum
(Theorem~\ref{limit-handles}) we must consider  
{\em handles of bounded thinness\/}, according to the terminology of
\cite{dmgm}; cf. Definitions~\ref{single-handle},
\ref{attachinghandle} and Remark~\ref{handlesdmgm}. This assumption is
necessary in that if we allow ``long-thin'' handles such as $
[-L , L]\times S^{d-1}(\e) $ (isometrically imbedded in $M$, as $\e
\downarrow 0$) then C.~Ann\'e \cite{anne1} proved that the spectrum
of the resulting manifold ``converges'' to the spectrum of $M$ {\em
plus\/} the spectrum of $d^2/dt^2$ on $ ( - L , L) $ with Dirichlet
boundary condition at $ t = \pm L$.

With different methods (and treating the case of smooth manifolds),
I. Chavel and E.A. Feldman considered in \cite{chavelfel4}, 
among other issues, the problem of attaching one handle to $M$ and
determining the size of this handle in order to have negligible
perturbations of the eigenvalues of $M$ (\cite[Theorem
5]{chavelfel4}). They, too, have to rule out long-thin handles and 
their condition involves the isoperimentric constant of the resulting
manifold. 

\vspace{.125in}
The organization of the paper is as follows. In the first section we
introduce the basic definitions and notation, define relaxed Lipschitz
manifolds (Definition~\ref{relax-mfld}) and make precise what we mean by
``representing a Lipschitz manifold in term of a relaxed Lipschitz
manifold'' (Definition~\ref{relaxriem}). 

In section~\ref{one-handle} we introduce and attach handles to $(S^d ,
g)$ (Definition~\ref{attachinghandle}). Then we represent the
resulting manifold by means of a relaxed Lipschitz 
manifold-with-boundary of simpler topological type
(Proposition~\ref{representation}).

Section~\ref{handle-repr} is devoted to the statements of our main
result, the proof of which 
is carried out in
Section~\ref{proof-of-the-theorem} and rests on Theorem~\ref{var-comp}
and Theorem~\ref{handles-comp}; because of the general nature of these
two results, we state and prove them separately in the Appendix,
\S~\ref{varcompmflds}.  

\vspace{.125in}
{\bf Acknowledgments.} It is a pleasure to thank Robert Gulliver for
many helpful suggestions and discussions during the preparation of the
paper and to acknoledge the hospitality of the School of
Mathematics at the University of Minnesota where the paper
was written. 

Thanks are also due to G~.Dal~Maso and U.~Mosco for their generosity
in providing the author with a copy of their preprint and for some
discussions. 

Financial support from the Italian Consiglio Nazionale delle Ricerche
during the preparation of the paper is also gratefully acknowledged.

\section{Notation \& Preliminaries.}
\noindent
When we define a quantity $A$ in term of other known quantities $B$,
we  write $A := B$.  

The {\em disjoint\/} union of two sets $A$, $B$ is denoted by $A \vee
B$. 

We consider real-valued functions;
if \( u \) is a function defined on some set \( X\), then \( u|_A \)
denotes its restriction to \( A \subset X\), 
\( u|_A : A \rightarrow \r\).  

We let $1_A(\cdot)$ denote the characteristic function of the set $A$,
namely $1_A(x) = 1$, if $x \in  A$, and zero otherwise.   

A measure \( \mu\) on a measurable space \( (X,{\cal A}) \) is a
countably additive function \( \mu: {\cal A} \rightarrow [0,+\infty]
\) defined on the \( \sigma\)-algebra \( {\cal A} \) which vanishes on
the empty set. Notice that we allow $\mu$ to assume the value \(
+\infty\).  If \( X\) is a topological space, then a Borel measure is
a measure defined on the \( \sigma\)-algebra of the Borel sets \( {\cal
B}\) of \(X\). We shall only consider Borel measure in the following
and sometimes for the sake of shortness we shall call them  measures.  

If \( X\) is a topological space, then the closure of a set $A \subset
X$ is denoted by $\overline A$; moreover \( A\subset \subset X\) means
\( A \subset \overline A \subset X\).  

\vspace{.125in}
Let $(X,d)$ be a  (separable) metric space, and let $B(x,\eta)$ denote
the open ball of center $x\in X$ and radius $\eta >0$.  

\begin{definition}[\( \eta \)-packing]\label{packings} Let $\eta > 0$;
a {\it  $\eta$-packing \/} in $(X,d)$ is  a collection of points $\{
x_i \}_{i\in I}$ in $(X,d)$ such that: 
\begin{itemize}
\item[(p$_1$)] $B(x_i,\eta) \cap B(x_j,\eta) = \emptyset$, for every
$i\not= j$, $i,j \in I$; 
\item[(p$_2$)] $\bigcup_{i\in I} B(x_i, 2\eta) = X$.
\end{itemize}
\end{definition}

\noindent
We can easily construct an $\eta$-packing in every metric space
$(X,d)$, in particular in a  manifold. Indeed, let $\eta > 0
$ be given, and let $x_1$ be any point in $X$; if the ball
$B(x_1,2\eta)$ covers $X$, then we are done. If not, there exists
$x_2$ whose distance from $x_1$ is greater than or equal to
$2\eta$. Therefore $B(x_1,\eta) \cap  B(x_2,\eta) = \emptyset$; if
$B(x_1,2\eta) \cup  B(x_2,2\eta) = X$, then we are done; otherwise
there is a point, which we call $x_3$ whose distance from both $x_1$
and $x_2$ is greater than or equal to $2\eta$, hence $B(x_i,\eta) \cap
B(x_3,\eta) = \emptyset$, $i=1,2$. Using  the induction  we find that
there exists an (at most countable) family $\{ x_i \}_{i\in I}$ which
satisfies the properties (p$_1$) and (p$_2$) above. Also, if  $X$ is
compact, then eventually a finite number of $x_i$'s will cover $X$. 

\begin{definition}
In a metric space  \( (X,d)\)  let  \( F_h: X
\longrightarrow [-\infty,+\infty] \) be a functional defined on \( X
\), \( h\in \n \). We say that the sequence \( (F_h) \) \(\Gamma\)-converges 
in  \( X 
\) to a functional \( F : X \longrightarrow [-\infty, +\infty] \) if
and only if
\begin{itemize}
\item[(a)] for every sequence \( (x_h) \) in \( X \), converging to \( x\in
X\) we have
\[
F(x) \le \liminf_{h\uparrow +\infty} F_h(x_h);
\]
\item[(b)] for every \( x \in X\) there exists a sequence \(
(\overline x_h) \) converging to \( x\) such that 
\[
\limsup_{h\uparrow+\infty} F_h(\overline x_h) \le F(x).
\]
\end{itemize}
We refer to Dal Maso's monograph \cite{dm1} for more information about \(
\Gamma\)-convergence.
\end{definition}

\begin{definition}\label{uniform-embedding}
Let $({\cal H}_h)_h$ be a sequence of Hilbert spaces; we say that
$({\cal H}_h)_h$ is {\em uniformly embedded\/} in the Hilbert
space ${\cal H}$ if there exists a linear, injective map
\[
I_h : {\cal H}_h \longrightarrow {\cal H}, \ \ h\in N
\]
such that the metric $\|\cdot \|_{{\cal H}}$ (induced by the Hilbert
structure of ${\cal H}$) 
on $ I_h ({\cal H}_h)$ is {\em uniformly\/} equivalent to the metric $
\| \cdot \|_{{\cal H}_h} $
of ${\cal H}_h$, i.e., there exists a constant $c_o$, possibly
depending on ${\cal H}$, but not on $h \in \n $,  such that 
\[
c_o^{-1} \| x - y \|_{{\cal H}_h} \le \| I_h(x) - I_h(y) \|_{{\cal H}}
\le c_o \| x - y \|_{{\cal H}_h}.
\]
\end{definition}

\vspace{.0625in}
In Theorem~\ref{limit-handles} we shall need the following definition,
which is a slight generalization of the convergence of subsets in a
Hilbert space as given in \cite[Definition 2.7.2]{mosco1}. 

\begin{definition}\label{conv-spaces}
Let $( {\cal H}_h)_h$ be a sequence of Hilbert spaces uniformly
embedded in ${\cal H}$ and let $S_h$ be a subset of ${\cal H}_h$,
$h\in \n$. 
We say that the sequence $(S_h)_h$ converges to a subset $S$
of ${\cal H}$ if the following conditions are satisfied: 
\begin{itemize}
\item[(i)] for every subsequence $(S_{h'})_{h'}$ of $(S_h)_h$ and for
every $x_{h'} \in I_{h'}(S_{h'})$ converging {\it  weakly\/} 
(in ${\cal H}$) to some $x\in {\cal H}$, we have $x \in S$; 
\item[(ii)] for every $x\in S$ there exists a sequence $(x_h)_h$ 
converging {\it  strongly\/} to $x$ such that $x_h \in I_h(S_h)$, for 
every $h \in \n$. 
\end{itemize}
\end{definition}

\vspace{.0625in}

In the following we let \( (M,g) \) denote a Lipschitz \(
d\)-dimensional manifold (with or 
without boundary), \( d\ge 2\), with \( g = (g_{ij})_{i,j=1}^d \) 
the metric tensor of \( M\).
We denote by \( \mbox{Vol}_g (dx) \) the {\em canonical measure\/} of
\( g\) on \( M\)  
\cite[Definition 3.90]{ghl}; in local coordinates \( x =
(x_1,\ldots,x_d) \) it can be written as  
\[ 
\mbox{Vol}_g(dx) = \sqrt{\det g\,}(x)~dx,
\]
where \( dx \) is the Lebesgue measure on \( {\bf R}^d\), so that 
\[
\mbox{Vol}_g(B) = \int_B \sqrt{\det g\,}(x)~dx,
\]
for every Borel set \( B\subset M\); the
term \( \sqrt{\det g\,}(x) \) sometimes will be referred to as the
local density of \( \mbox{Vol}_g\). 

We occasionally will also use the notation $ds^2_M$, $\mbox{Vol}_M$ to
denote respectively the metric of $M$ and its canonical measure. 

\vspace{.0625in}
We point out that Definition 3.90 in \cite{ghl} requires the manifold
to be $C^\infty$; however their definition can be extended in our
Lipschitz framework. 

\vspace{.0625in}
We  define \( L^2(M,g) \) as the space of all functions on $M$ such
that  
\[
 \int_M f^2\mbox{Vol}_g(dx) < +\infty.
\]

\vspace{.0625in}\noindent
In the framework of Lipschitz manifolds it still makes sense to consider
the {\em Dirichlet functional\/}  on $ (M,g)$ 
(\cite{eisenhart}) 
$$
{\cal D}_{M,g}(u)  := \int_M  \sum_{i,j=1}^d g^{ij}D_iu D_j u
\mbox{Vol}_g(dx), 
$$
for Lipschitz  functions $u$ defined on $M$; notice that $D_i u :=
\partial u/\partial 
x_i$, $i=1,\ldots, d$, is defined $\mbox{Vol}_g$-almost everywhere.

\begin{definition}\label{hone}
We let $H^1(M,g)$ denote the closure of $\mbox{Lip}(M)$, the family of
Lipschitz functions on $M$, under the norm induced by  
$$
{\cal D}_{M,g} (u) + \int_M u^2 \mbox{Vol}_g(dx).
$$
\end{definition}

\vspace{.0625in}\noindent
For $\lambda >0$, and $f\in L^2(M,g)$, we can introduce the functional
$F^\lambda : L^2(M,g) \longrightarrow [0,+\infty]$ defined by 
$$
F^\lambda(u) : = {\cal D}_{M,g} (u) + \lambda \int_M u^2 \mbox{Vol}_g(dx) -
2\int_M f u  \mbox{Vol}_g(dx),
$$
if $u\in H^1(M,g)$ and $F^\lambda(u) := +\infty$ otherwise in
$L^2(M,g)$. This functional is strictly convex, coercive and lower
semi-continuous in the 
strong topology of $L^2(M,g)$ hence, by the Direct Method in the
Calculus of Variation, it has a unique minimum point $u_f$. 

\vspace{.125in}
Via standard arguments it is possible to prove the following result. 

\begin{prop}\label{resolvent}
The resolvent operator $R^\lambda : L^2(M,g) \longrightarrow
L^2(M,g)$, which 
associates to every $f\in L^2(M,g)$ the minimum
point $u_f$, is a compact, positive, self-adjoint operator.
\end{prop}

Thus the resolvent operator has a sequence of proper values
$(\sigma^o_i)_{i\in \n}$ having zero as accumulation point. 

\begin{definition}\label{hone-spectrum}
The spectrum of $(M,g)$ is the sequence of the proper values 
$(\sigma^o_i)_{i\in \n}$ of the resolvent operator $R^\lambda$.
\end{definition}

If  $( M, g)$ is smooth, then let $-\Delta_g$ denote its
Laplace-Beltrami operator.  

\begin{definition}\label{nat-eigen}
We say that $\alpha$ is an eigenvalue of $-\Delta_g$ if
\[
\begin{cases} 
-\Delta_g u = \alpha u,\ \  \text{on $M$}\\
\\
\mbox{natural boundary condition}
\end{cases}
\]
where we agree that ``natural boundary condition'' means
homogeneous Neumann boundary condition if $\partial M \not=
\emptyset$, and no boundary condition if $\partial M =\emptyset$.
\end{definition}


\begin{remark}\label{remark-hone}
Notice that if $(M,g)$ is a smooth manifold,
then $\sigma^o_i = (\lambda + \lambda_i)^{-1}$, where
$(\lambda_i)_{i\in \n}$ is the sequence of eigenvalues of
$-\Delta_g$.   
\end{remark}


\begin{definition}\label{lscr}
In what follows we  shall need the {\em lower semi-continuous
regularization\/} in $L^2(M, \mbox{Vol}_g)$ of the functional \( {\cal
D}_{M,g}(\cdot) \):
\begin{equation}
\widetilde{{\cal D}}_{M,g}(u) := \begin{cases}
{\cal D}_{M,g} (u), \ \  \mbox{if} \ \ u\in H^1(M,g) \\   
\\ 
+\infty, \ \ \mbox{otherwise in}\ L^2(M,g).
\end{cases} 
\label{lsc}
\end{equation}
\end{definition}

\begin{definition}\label{mg-capacity}
Let \( A \) be a given relatively compact open set contained in
\( \bigl( M,g \bigr) \), and let \( E \) be a Borel set, \( E\subset
A\). Then the  \( (M,g)\)-capacity of \( E\) with respect to \( A\)
is defined by 
\[
(M,g)\mbox{-cap}(E,A) := \inf \left\{ {\cal D}_{M,g}(u) : u\in 
{\cal K}(E,A) \right\},
\]
where \( {\cal K}(E,A) := \{ u\in\mbox{Lip}_o(A): u \ge 1\ \mbox{on a 
neighborhood of}\ E\} \), and $\mbox{Lip}_o(A)$ denotes the family of
all Lipschitz functions whose support is contained in $A$.  
\end{definition}

\begin{remark}\label{cap-zero}
Notice that the property of having capacity zero is unchanged if \(
A\) is replaced by a larger relatively compact open subset of \(
M\). Therefore we shall say that a property \( P(x) \) holds \(
(M,g)\)-{\em quasi everywhere\/} (q.e.), or for \( (M,g)\)-{\em quasi
every\/} \( x\in M\), if the set 
\[
\{ x\in M : P(x) \ \mbox{is not true}\ \}
\]
has \( (M,g)\)-capacity zero. 
\end{remark}


\begin{remark}\label{quasi-repr}
Modifying suitably some arguments in \cite{federerziemer}, it is
possible to show that for every \( u\in H^1(M,g) \)
there exists \( \widetilde u: M \rightarrow \r \) such that
\[
\lim_{r\downarrow 0} \frac{1}{\mbox{Vol}_g(B(x,r))} \int_{B(x,r)}
|\widetilde u(x) - u(y)|~\mbox{Vol}_g(dy) = 0,
\]
and \( \widetilde u = u \) \( \mbox{Vol}_g\)-almost everywhere on \(
M\). The function $\tilde u$ is uniquely determined up to a set of
$(M,g)$-capacity zero and is continuous when restricted to the
complement of open sets with arbitrarily small $(M,g)$-capacity. In
the following we shall identify each \( u\in H^1(M,g) \) with \(
\widetilde u \), so that we can say that \( u \) is determined \(
(M,g) \)-quasi everywhere.
\end{remark}


\begin{definition}\label{mzero}
Let \( {\cal M}_o(M,g) \) denote the class of all Borel measures \(
\mu \) on \( M \) which are {\em absolutely continuous\/} w.r.t. the
\( (M,g)\)-capacity, \ie  
\[
\mu(E) = 0 \ \ \mbox{whenever}\ E\ \mbox{has \( (M,g)\)-capacity
zero}. 
\]
We introduce an equivalence relation among measures in \( {\cal
M}_o(M, g)\): We say that \( \mu, \nu\in {\cal M}_o(M,g) \) are
{\em equivalent\/}, \( \mu \sim \nu\), if and only if 
\[
\int_M u^2 d\mu = \int_M u^2 d\nu,
\]
for every \( u\in H^1(M,g)\). By Remark~\ref{quasi-repr} above each
function $u \in H^1(M,g)$ is determined $(M,g)$-quasi-everywhere,
hence  \( \mu\)- and \( \nu\)-almost everywhere. Thus the integrals
above are well-defined, possibly equal to \( +\infty \).
\end{definition}

\begin{definition}
\label{inftymeasure}
Let \( E \) be a
given Borel set; then 
\[
\infty_E(B) := \left\{
\begin{array}{ll}
0, \ \mbox{if}\ B\cap E \ \mbox{has \( (M,g)\)-capacity zero;} \\
\\
+\infty, \ \mbox{otherwise.} 
\end{array}\right.
\]
Thus \( \infty_E(\cdot) \in {\cal M}_o(M,g)\). A function \( v \in
H^1(M,g) \) belongs to \( L^2(M,\infty_E(dx)) \) if and only if \( v =
0 \) \( (M,g)\)-quasi everywhere on \( E\).
\end{definition}
\begin{remark}\label{mz-vol}
The canonical measure \( \mbox{Vol}_g\) belongs to \( {\cal M}_o(M,g)
\). Indeed, let \( E\) be a Borel set which has capacity zero; then \(
(M,g)\mbox{-cap}(E,A) = 0\), for some relatively compact open set \(
A\subset M\). Therefore, by definition, for every \( \e> 0\) there is \(
\varphi_\e \in \mbox{Lip}(A)\), with the support contained in $A$,
such that \( {\cal D}_M(\varphi_\e) < \e\), and \( \varphi_\e \ge 1\)
on a neighborhood of \( E \). Thus the characteristic function of \( 
E\) is less than or equal to \( \bigl(\varphi_\e\bigr)^2 \) on \( A
\), hence \( \mbox{Vol}_g(E) \le \int_A \bigl(\varphi_\e(x)\bigr)^2
\mbox{Vol}_g(dx) \). By the Poincar\'e inequality
\[
\int_A \bigl(\varphi_\e(x)\bigr)^2 \mbox{Vol}_g(dx) \le C {\cal
D}_M(\varphi_\e),
\]
where the constant \( C\) possibly depends on \(A\), we get 
\[
\mbox{Vol}_g(E) \le {\cal D}_M(\varphi_\e) < \e,
\]
which implies \( \mbox{Vol}_g(E) = 0\), by the arbitrariness of \(
\e>0 \). 
\end{remark}

\begin{definition}[\cite{dmgm}]\label{relax-mfld}
A {\em relaxed\/} Lipschitz manifold is by definition a
4-tuple \( \bigl(M,ds^2_M,T,\mu \bigr) \) where \( (M,ds^2_M) \) is a
Lipschitz manifold, possibly with boundary; \( T: M \rightarrow M \)
is an isometry, \( 
T\circ T = \mbox{id}_M\) and the fixed-point set 
of \(T\), \( \mbox{Fix}(T) = \{ x\in M : T(x) = x\} \), is a
submanifold of \( (M,ds_M^2) \); \( \mu \) is a measure which belongs
to \( {\cal M}_o(M,g) \). 
\end{definition}

\vspace{.0625in}\noindent
We have the following general definition (cf. Definition 1.8 in
\cite{dmgm}).

\begin{definition}\label{relaxriem}
Let \( (\overline M,ds_M^2) \) be a manifold-with-boundary
and  \( (N,ds_N^2) \) be a manifold (possibly with \(
\partial N = \emptyset\)); let \( \mbox{Vol}_M \), \( \mbox{Vol}_N
\) be the canonical measures, and  \( {\cal D}_M(\cdot) \), \( {\cal
D}_N(\cdot) \) be the Dirichlet functionals of respectively \( M \)
and  \( N \). Let \( T: \overline M \longrightarrow \overline M \) be
an  isometry such that \( T\circ T = \mbox{id}_{\overline M} \), and
the   fixed-point set 
\( \mbox{Fix}(T) \) is a submanifold of \( M \) (possibly \(
\mbox{Fix}(T) = \emptyset \)). Furthermore, let \( \nu \) 
be a Borel measure on \( \overline M \). We say that the manifold \(
(N,ds_N^2) 
\) is {\em represented\/} by the relaxed manifold \(
\bigl(\overline M, ds_M^2, T , \nu\bigr) \) if there is an isometry
\[ 
I : L^2(N,\mbox{Vol}_N) \longrightarrow L^2(M,\mbox{Vol}_M) 
\]
and for \( v\in H^1(N) \) we have 
\[
{\cal D}_N(v) = {\cal D}_M(u) + \int_{\overline M} \bigl[u(x) -
u(T(x))\bigr]^2 \nu(dx),
\]
with \( u = I(v) \).
\end{definition}

\begin{remark}\label{finiteness}
Notice that \( {\cal D}_Y(v) < +\infty\) if and
only if \( {\cal D}_X(u) <+\infty \) and \( u(\cdot) - u( T(\cdot))
\in L^2(X, \nu(dx)) \). 
\end{remark}

\section{Handles}\label{one-handle}

In this section we describe an explicit procedure to attach
a handle to $S^d$, and denote the resulting manifold by $(N, g^N)$;
the topological type of $N$, then,  will be that of $S^{d-1} \times
S^1$. Then we show that $(N , g^N)$ can be represented (in the sense
of Definition~\ref{relaxriem} above) by a relaxed manifold $(M_1, g_1,
T, \infty_{\partial M_1})$, where $M_1$ is homeomorphic to $S^d$ wth two
punctures. 

\vspace{.0625in}
We point out that the handles we shall consider here
are called ``handles of bounded thinness'' in \cite[\S 3]{dmgm}.

\vspace{.125in}\noindent
As a matter of notation let \(
S^{d-1} \) be the \( (d-1) \)-sphere of radius 1, and 
\( S^{d-1}(t) \) be the \( d-1 \)-sphere of radius \( t > 0  \). 

Let $0 < \e <1$; let us consider the cylinder 
\[
H := [-1,1] \times S^{d-1}(\e),
\]
with its usual product metric 
\begin{equation}
ds_H^2 = dy^2 + \e^2 ds^2_{S^{d-1}}(\omega).
\label{metrich}
\end{equation} 

\noindent
Let us consider the map
\begin{equation}
\begin{array}{cc}
T: {\r}^{d+1} \longrightarrow {\r}^{d+1} \\
T(x) : = - x
\end{array}
\label{antipodal-map}
\end{equation}
Notice that \( T \) is an isometry, \(T\circ T = \mbox{id}_{{\r}^{d+1}}
\), and the fixed-point set of \( T\) 
reduces to the singleton \( \{ 0\} \). In the following we are only
concerned with the restriction of this map to \( S^d \), and 
in order not to have too heavy notation we still denote this
restriction by
\( T\). We notice that \( T\) (the perhaps familiar 
antipodal map on \( S^d \)) maps \( S^d \) onto itself, is an
isometry, \( T\circ T = \mbox{id}_{S^d}\) and  the fixed-point set
is empty. 

\vspace{.0625in}
We shall find useful in the following to employ
cylindrical coordinates so that if \( x\in {\r}^{d+1}
\), then \( x = (y,r\omega) \) where \( y\in \r\), \( \omega \in
S^{d-1} \), \( r\ge 0 \), and \( |x|^2 = y^2 + r^2\); in
particular, 
\[
S^d = \{ (y,r\omega)\in{\r}^{d+1} : y^2 + r^2 =1 \},
\]
and the metric of \( S^d \) is in these coordinates 
\begin{equation}
ds_{S^d}^2 = (1-r^2)^{-1} d r^2 + r^2 ds^2_{S^{d-1}}(\omega).
\label{metrics2} 
\end{equation}
Also, \( H = \{ (y,r\omega) : y\in [-1,1], r = \e\} \). We shall also
consider the map  
\begin{equation} 
\begin{array}{cc}
T_H : H \longrightarrow H \\
(y,\e\omega) \mapsto (-y,\e\omega)
\end{array}
\label{handle-map}
\end{equation}
Notice that \( T_H\) is an isometry, \( T_H \circ T_H = \mbox{id}_H\),
\( \mbox{Fix}(T_h) = \{ 0 \} \times S^{d-1}(\e) \).

\begin{definition}\label{lipschitzfamily}
Let \( \Theta \in [1,+\infty)\). Let \( {\cal F}_\Theta \) 
be the family of all pairs \( (\e, r) \) where \( \e \in (0,1) \), and
\( r(\cdot) \) is an increasing, Lipschitz
continuous function, \( r : [0,1] \longrightarrow (0,1) \), 
with 
\[
0 < D^{-1} < \frac{dr}{dt}(t) < D, \ \mbox{a.e.\ } t\in (0,1),
\]
such that
\[
\max\left\{ \frac{\e}{r(0)}, \frac{r(1)}{\e},
\frac{D}{\sqrt{1-\bigl(r(1)\bigr)^2\,}}\right\} < \Theta.
\]
\end{definition}

\begin{definition}[Handle]\label{single-handle}
Let \( (\e, r_1)\in{\cal F}_\Theta \), and let 
$\varphi : H \setminus \mbox{Fix}(T_H)
\longrightarrow S^d$ be the function defined by 
\begin{equation}
\varphi : (y, \e\omega) \mapsto 
(\mbox{sgn}(y)\sqrt{1 - (r_1(|y|))^2}, 
\mbox{sgn}(y)r_1(|y|)\omega). 
\label{handlephi}
\end{equation}
The pair \( (H,\varphi) \) is by definition a {\em handle\/}.  
\end{definition}

\begin{ex}\label{handle-line}
Let \( \Theta >1\) be given. Let \( \e > 0 \), and let us consider 
\[
r_1(t) := \beta t + \alpha(1-t), \ \ t\in [0,1].
\]
We choose \( \alpha = \e\delta_o \), \( \beta = \e\delta_1 \) with 
\( 0 < \delta_o < \delta_1 < 1\); if \( \delta_o > 1/\Theta \), then
\( (\e, r_1)\in {\cal F}_\Theta\), hence the pair \( (H,\varphi) 
\), where \( \varphi \) is as in Definition~\ref{single-handle} above,
is a handle. 
\end{ex}

\noindent
Let 
\begin{align}
U &:= B((1,0), \sin^{-1}r_1(1)) \vee
B((-1,0), \sin^{-1}r_1(1)), \label{holes-annuli1}
\\
E &:=  \overline{ B((1,0), \sin^{-1}r_1(0))} \vee
\overline{ B((-1,0), \sin^{-1}r_1(0))}. \label{holes-annuli2}
\end{align}
Notice that  
\begin{eqnarray*}
\varphi(\mbox{int}(H) \setminus \mbox{Fix}(T_H))  &\!\!\! =\!\!\! & 
\varphi( \{ (y,\e\omega)\in H : y \not= 0\}) \subset U, \\ 
\varphi(\partial H) & \!\!\! = \!\!\! & \varphi( \{ \pm 1\} \times 
S^{d-1}(\e)) =  \partial U, \\
E & \!\!\! = \!\!\!& U\setminus \varphi(\mbox{int}(H)\setminus
\mbox{Fix}(T_H)) 
\end{eqnarray*}

\begin{definition}[Attaching the handle]\label{attachinghandle}
To attach the handle $(H,\varphi)$ to
$S^d$ we consider \( (S^d \setminus E) 
\cup H \) and {\em identify\/} each point \( z\in H\setminus
\mbox{Fix}(T_H) \) with \( \varphi(z) \in U\setminus E \): \( z \sim
\varphi(z) \). Then the {\em topological type\/} of 
\[
N := \bigl( (S^d \setminus E) \cup H\bigr)/\sim
\]
is that of $S^{d-1} \times S^1(1)$. The {\em Lipschitz metric\/} of
the resulting manifold $N$ is by definition the metric of 
\( S^d \) on \( N\setminus E \) and that of \( H \) on \(
\mbox{int}(H) \): 
\begin{equation}
ds_N^2 := 
\left\{\begin{array}{ll}
ds_{S^d}^2, & \ \mbox{on}\ S^d \setminus U, \\
\\
ds_H^2, & \ \mbox{on}\ \mbox{int}(H).
\end{array}\right.
\label{resultingmetric}
\end{equation}
\end{definition}

\begin{remark}\label{boundary}
Notice that the boundary of \( H \), \( \partial H = \{ \pm 1\} \times
S^{d-1}(\e) \), is thus identified with the hypersurfaces \( \{ (y,
r\omega)\in S^d : r = r_1(1) \} \).
\end{remark}

\begin{remark}\label{handlesdmgm}
It can be shown (cf. Remark 3.6 in \cite{dmgm}) that the metric \(
g^N = (g^N_{ij})_{i,j} \) on the resulting manifold \( N \) is 
{\em uniformly bounded\/} by the metric \(
( g_{ij})_{i,j} \) of \( S^d \), {\em i.e.\/}, 
\[
\Theta^{-2} \sum_{i,j} g_{ij}\xi^i\xi^j \le
\sum_{i,j}g^N_{ij}\xi^i\xi^j 
\le \Theta^2 \sum_{i,j} g_{ij}\xi^i\xi^j, 
\]
for every \( \xi\in{\r^n} \). 

Moreover the complement of \( E\) is {\em uniformly strongly 
connected\/}, {\em i.e.\/}, for each $u\in H^1(S^d \setminus E,g)$
there is an extension of \( u \), $v = \pi u$, with $v\in
H^1(S^d,g)$, and  
\[
\| v \|_{H^1(S^d,g)} \le \Theta^{1/2} \| u \|_{H^1(S^d 
\setminus E,g)}.
\]
The constant \( \Theta \) is the same that appears in
Definition~\ref{single-handle}.
\end{remark}

\begin{remark}\label{handle-line-ex}
Let \( (H,\varphi) \) be the handle with the function \( r_1(\cdot) \)
as in Example~\ref{handle-line}. Then comparing (\ref{metrics2}),
(\ref{metrich}), (\ref{resultingmetric}), it 
is not difficult to see that the metric on \( N \) is continuous 
({\em i.e.\/}, the manifold is \( C^1 \)) if and only if \( B = r_1(1)
= \e\). We can also have the resulting manifold to be of class \(
C^{k+1} \) provided that we modify the metric of the handle near \(
S^{d-1}(\e) \times \{ \pm 1\} \) so that the new metric is of class \(
C^k \) and as \( t\to \pm 1\) the metric is \( ds^2_H = \bigl( 1-
r_1(t) \bigr)^{-1/2}\bigl( dr_1/dt\bigr)^2 dt^2 + \bigl(
r_1(t)\bigr)^2 ds^2_{S^{d-1}}\) plus terms which vanish at \( t=
\pm 1\) together with all derivatives of order less than or equal to
\(k\). 

If instead \( r_1(1) \not= \e \), then the metric is only
piecewise continuous: In fact, the metric has a
discontinuity along the hypersurfaces \( \{ (y,r\omega) \in S^d :  r
= r_1(1), y = \pm\sqrt{1 - r^2\,}\} \) 
where \( \partial H \) is attached to \( S^d \), whereas the metric
is smooth elsewhere on \( S^d \). Therefore, if \(
r_1(1)  \not= \e \), \( N \) is a Lipschitz, but not a \( C^1 \), manifold. 

Notice that the sectional curvature may present a {\em
distributional component\/} along the hypersurfaces 
\( \{ (y,r\omega) \in S^d :  r = r_1(1), y = \pm\sqrt{1 - r^2\,}
\). 
\end{remark}

\noindent
Let \( \overline M_1 \) be the manifold-with-boundary \( S^d
\setminus E \), \( \partial M_1 = \partial E \), whose metric \(
g_1 = (g_{ij,1})_{i,j} \) is equal to
\begin{equation}
(g_{ij,1})_{i,j} := 
\left\{\begin{array}{ll}
(g_{ij})_{i,j}, & \ \mbox{on}\ S^d \setminus U, \\
\\
( g^H_{ij})_{i,j}, & \ \mbox{on}\ U\setminus E,
\end{array}\right. 
\label{relaxmetric}
\end{equation}
where \( (g_{ij})_{i,j} \) denotes the metric of \( S^d \), 
and 
\( ({g}^H_{ij})_{i,j} \) is the pull-back of the metric of
\( H \) under the map \( \varphi^{-1} :
\varphi(H \setminus \mbox{Fix}(T_H)) \longrightarrow H \). 
We notice that the metric in (\ref{relaxmetric}) is {\em uniformly
bounded\/} by the metric \( g = (g_{ij})_{i,j} \) of \( S^d \). 
We point out that the map \( T \) introduced in (\ref{antipodal-map})
maps \( \overline M_1\) onto itself, is an isometry there and \(
\mbox{Fix}(T) = \emptyset \).

\vspace{.125in}\noindent
We have the following result (cf. Lemma~1.11 in \cite{dmgm}). 

\begin{prop}\label{representation} 
The manifold \( ( N, g^N ) \) is
represented (in the sense of Definition~\ref{relaxriem}) by the
relaxed manifold-with-boundary \( \bigl( \overline M_1,
g_1, T, \infty_{\partial M_1}  \bigr) \). 
\end{prop}

\noindent
In particular, for each \( v\in H^1(N,g^N) \) there exists \( u \in
H^1(\Omega,a)\) such that
\begin{equation}
{\cal D}_N(v) = {\cal D}_{M_1}(u) + 
\int_{\overline M_1} \bigl[ u(x) - u(T(x))
\bigr]^2 \infty_{\partial M_1}(dx).
\label{dirichletrep} 
\end{equation}

\begin{remark}\label{cutandpaste}
Notice that the right hand side in (\ref{dirichletrep}) is finite if 
and only if 
\begin{equation}
u(\cdot) = u(T(\cdot)), \ \mbox{q.e. on~~} \partial M_1. 
\label{glueing}
\end{equation}
(cf. Definition~\ref{inftymeasure}.)
Roughly speaking, we may say that the
representation of \( N \) has been obtained by the following 
``cut-and-paste'' procedure: We cut the handle of $N$  along
$\mbox{Fix}(T)$ and get the manifold-with-boundary $\overline M_1$;
notice that $\overline M_1$ is homeomorphic to $S^d$ with two
punctures.
The presence of the handle, as far as the Dirichlet
functional is concerned, is then represented by the non-local term 
\[
\int_{\overline M_1} \bigl[ u(x) - u(T(x))
\bigr]^2 \infty_{\partial M_1}(dx),
\]
which, via (\ref{glueing}), glues together the two components of \(
\partial M_1 = \partial E \) and gives back the handle.

If the measure $\mu$ appearing in the non-local term above is 
finite, with $(N, g^N)$ represented by $(\overline M_1, g_1, T, \mu)$,
then by analogy with the case of the ``infinite'' measure, we may say
that the handle is ``weakly'' (or ``partially'') attached to $N$.
\end{remark}


\section{The main results} \label{handle-repr}

\noindent
Let $(\eta_h)_h$ be a sequence of positive numbers, $\lim_{h\uparrow +
\infty} \eta_h = 0$ and  let $(x_i)_{i\in I_h}$ be an $\eta_h$-package
in $(S^d,g)$ (cf. Definition~\ref{packings}),  for each  $h\in\n$.
Let  also  $\overline x_i := T(x_i)$, $i\in I_h$, where $T$ is the
antipodal map on $S^d$.  

Let \( (\e_h)_h \) be a sequence of positive numbers such that 
\begin{equation}
\lim_{h\uparrow+\infty} \eta_h / \e_h  = 0.
\label{pinco}
\end{equation}
Let us consider  \( H_h : = [-1,1]\times S^{d-1} ( {\e_h} ) \), whose
metric is equal to $ds^2_{H_h} = \e^2_h\left( dy^2 + ds^2_{S^{d-1}}
\right)$. 

Using cylindrical coordinates we can assume
that  \( x_i = (1,0) \) and \(\overline x_i = T(x_i) = (-1,0) \), for
$i\in I_h$. Hence,  for a given $\Theta > 1$, we can define  maps
(similarly as in the previous 
\S~\ref{one-handle}) 
\[
\begin{array}{ll}
\varphi^\pm_i : H_h\cap \{ \pm y > 0 \} \longrightarrow  S^d \\
\\
\varphi^\pm_i (y,\e_h\omega) =  
( \pm \sqrt{1 - (r_{\e_h}(\pm y))^2}, 
\pm r_{\e_h} (\pm y)\omega),
\end{array}
\]
with (cf. Example~\ref{handle-line})
\[
r_{\e_h}(t) := (\delta_1\e_h\eta_h)t + (\delta_o\e_h\eta_h)(1-t), \ \ \ 
0 < \frac{1}{\Theta} < \delta_o < \delta_1 <1, \ \ \ t\in [0,1],
\]
so that  $ ( r_{ \e_h} , \e_h\eta_h )\in {\cal F}_\Theta $;  if we define 
\[
\varphi_{i,h}(\cdot) := \left\{
\begin{array}{ll}
\varphi^+_{i,h}(\cdot), \ \mbox{on}\ H_h\cap \{ y>0 \},\\
\\
\varphi^-_{i,h}(\cdot), \ \mbox{on}\ H_h\cap \{ y<0 \},
\end{array}\right.
\]
then the pair \( (H_h, \varphi_{i,h}) \) is a handle, as in
Definition~\ref{single-handle}. 

\begin{definition}
\label{h-attachinghandles} %
(cf. Definition~\ref{attachinghandle}) 
We let $(N_h,g^{N_h})$ denote 
the manifold obtained by attaching the handles
$(H_h,\varphi_{i,h}) $, $i\in I_h$, to $S^d$. 
\end{definition}

\begin{remark}\label{lip-mfld}
As $r_{\e_h}(1) = \eta_h\e_h \not= \e_h$, then the metric $ g^{N_h} $
is not continuous, but only piecewise continuous;
cf. Remark~\ref{handle-line-ex}. Thus $(N_h,g^{N_h})$ is a
Lipschitz manifold, but not a $C^1$ manifold.  
\end{remark}

\begin{definition}
\label{h-spectrum}
(cf. Definition~\ref{hone}) We let $(\sigma^h_i)_{i\in \n}$ denote the
spectrum of $(N_h,g^{N_h})$. 
\end{definition}

\vspace{.125in}\noindent
Our main result is the following.

\begin{th}\label{limit-handles}
Let $r_h : = \arcsin r_{\e_h} ( 0 )$, and  assume that  
\begin{equation}
\frac{\alpha}{2} := \begin{cases}
\displaystyle {
\lim_{h\uparrow + \infty}  \frac{r_h^{d-2}}{\eta_h^d}} &  \ \mbox{if}\ \
d\ge 3,
\\ 
\\
\displaystyle{
\lim_{h\uparrow + \infty} \frac{-1}{\eta^2_h\log r_h }} & \ \mbox{if}\ \
d=2,
\end{cases}\label{hypo-radii}
\end{equation}
with $ 0 \le \alpha <  +\infty$.  Let moreover $(\sigma_i)_i$ be the
sequence of the proper values  of the resolvent operator
$R^\lambda_\infty : L^2 (S^d , g ) \rightarrow L^2 ( S^d , g ) $
corresponding to the functional  
\begin{align*}
{\cal D}_\infty(u) := 
{\cal D}_{S^d} ( u )  & + \frac{ \alpha }{ 2^d }\int_{ S^d } \bigl[
u(x) - u(T(x)) \bigr]^2 \mbox{Vol}_g(dx) + \\ 
& + \lambda \int_{ S^d } u^2\mbox{Vol}_g(dx) - 2\int_{ S^d }
fu\mbox{Vol}_g(dx),  
\end{align*}
for $u\in H^1(S^d, g)$; then 
$$
\lim_{h\uparrow +\infty} \sigma_i^h = \sigma_i,\ \ \ \mbox{for each }\
i \in \n. 
$$
Moreover the sequence $(L^2(N_h ,g^{N_h} ))_h$ is uniformly embedded
in $L^2(S^d,g)$ and if $r \le i \le s$ are such that $
\sigma^h_{s-1} < \sigma^h_{ s }  =   \sigma^h_{i} =  
 \sigma^h_{ r }  <  \sigma^h_{ r +1}$, then the linear subspace
spanned by $ \{ u^h_r, \ldots,  u^h_s
\}$ in $L^2(N_h ,g^{N_h} )$ converges (in the sense of
Definition~\ref{conv-spaces}) to the eigenspace corresponding to
$\sigma_i$, where $u_r^h,\ldots, u^h_s$ are the eigenfunctions
corresponding to $\sigma^h_i$. 
\end{th}

\vspace{.125in}\noindent
Let  $u$ be a function defined on $S^d$;  following \cite[\S
5]{dmgm} let us define the respectively odd and even part of $u$ by  
\begin{align*}
u_{\rm odd}(\cdot) & := \frac{1}{2}[u(\cdot) - u(T(\cdot)) ] \\
u_{\rm even}(\cdot) & := \frac{1}{2}[u(\cdot) + u(T(\cdot)) ],
\end{align*}
so that $u(\cdot) = u_{\rm odd}(\cdot)  + u_{\rm even}(\cdot) $. 
If moreover $\Xi$ is a space of functions defined on $S^d$, then we
denote by $\Xi_{\rm odd}$ (resp. $\Xi_{\rm even}$) the subspace of
$\Xi$ consisting of all odd (resp. even) functions on $S^d$. 

It can be shown that both the standard metric $g$ and canonical measure
$\mbox{Vol}_g$ of $S^d$ are invariant under the action of the
antipodal map $T$. Thus it can be proved that
both $L^2(S^d,g)$ and $H^1(S^d,g)$ split into their even and odd
part as follows 
\begin{align}
L^2(S^d,g) & = L^2_{\rm odd}(S^d,g) \oplus L^2_{\rm
even}(S^d,g)  \label{orth-dec-1}\\
H^1(S^d,g) & = H^1_{\rm odd}(S^d,g) \oplus H^1_{\rm
even}(S^d,g). \label{orth-dec-2}
\end{align}
(The orthogonal decomposition is with respect to the inner product of
$L^2(S^d , g)$ in the former, and in $H^1(S^d,g)$ in the latter.)

Thus the limit functional ${\cal D}_\infty(\cdot)$ can
be written as  
$$
{\cal D}_\infty(u) = {\cal D}_{\infty, {\rm even }} (u_{\rm even})  +
{\cal D}_{\infty, {\rm odd }}(u_{\rm odd}), 
$$
where
\begin{align*}
 {\cal D}_{\infty, {\rm odd }} (u_{\rm odd}) & = {\cal D}_{S^d} (
u_{\rm odd} ) + \\ 
& \left(  \lambda + \frac{\alpha}{2^d} \right)
\int_{ S^d } u_{\rm odd}^2\mbox{Vol}_g(dx) - 2\int_{ S^d }
u_{\rm odd} f_{\rm odd}\mbox{Vol}_g(dx),
\end{align*}
and
\begin{align*}
{\cal D}_{\infty, {\rm even}} (u_{\rm even} &  = {\cal D}_{S^d} (
u_{\rm even} ) + \\ 
&  \lambda \int_{ S^d } u_{\rm even}^2\mbox{Vol}_g(dx) - 2\int_{ S^d }
u_{\rm even} f_{\rm even}\mbox{Vol}_g(dx) 
\end{align*}
Using (\ref{orth-dec-1}), (\ref{orth-dec-2}) above, we have that
the Euler equation associated  with  ${\cal D}_{\infty}(\cdot)$ can be
``decoupled''  into two equations, one for the odd part and the other
for the even part of ${\cal D}_\infty(\cdot)$ as follows:
\begin{equation}
\begin{cases}
\displaystyle{
- \Delta_g u_{\rm odd} + \left( \lambda + \frac{\alpha}{2^d} \right)
u_{\rm odd} = f_{\rm odd}, }  &  
\mbox{ in $S^d$}\\
u_{\rm odd} \in H^1(S^d,g),
\end{cases}\label{euler-odd}
\end{equation}
and 
\begin{equation}
\begin{cases}
\displaystyle{
- \Delta_g u_{\rm even} +  \lambda   u_{\rm even} = f_{\rm even}, }  & 
\mbox{ in $S^d$}\\
u_{\rm even} \in H^1(S^d,g).
\end{cases}\label{euler-even}
\end{equation}

\begin{remark}\label{odd-even}
The decoupling of the Euler equation associated with the limit
functional ${\cal D}_\infty (\cdot)$ into (\ref{euler-odd}) and
(\ref{euler-even})  implies that the sequence
$(\sigma_i)_{i\in \n} $ of proper values of the resolvent operator
$R^\lambda_\infty$ splits 
into two sequences $(\sigma^{\rm odd}_i)_{i\in \n} $,  $(\sigma^{\rm
even}_i)_{i\in \n} $: the odd and even part of the
spectrum. Thus we see that adding an
increasing number of handles affects only the odd part of the spectrum
with the occurence of the Lenz shift phenomenon in (\ref{euler-odd}).
\end{remark}

\section{Proof of the main result}\label{proof-of-the-theorem}
\noindent
The strategy to prove our main result, Theorem~\ref{limit-handles}, is
as follows: We introduce the relaxed manifold \( (\overline
M_h,g_h,T,\infty_{\partial M_h}) \) 
(Definition~\ref{h-relaxed-mfld}) and prove that \( (\overline
M_h,g_h,T,\infty_{\partial M_h}) \) represents $(N^h, g^{N_h})$
(Proposition~\ref{h-representation}). Thus the value of the
functional $F^\lambda_h ( \cdot )$ introduced above is equal to the
value of the functional ${\cal R}_h^\lambda : L^2 ( S^d , g )
\rightarrow [0, +\infty] $, introduced in
Definition~\ref{relaxed-funct}. Finally 
we use Theorem~\ref{var-comp}, Theorem~\ref{handles-comp} and an
adaptation to our framework of a derivation-type argument in
\cite{bdmm2} to conclude the proof of Theorem~\ref{limit-handles}. 

\vspace{.125in}\noindent
Let $r_h = \arcsin r_{\e_h}(0)$, $R_h =  \arcsin r_{\e_h}(1)$
(cf. (\ref{holes-annuli1}), (\ref{holes-annuli2})); define
\begin{align}
E_h & :=   \bigcup_{i\in I_h}\left[
\overline{B(x_i,r_h)} \vee \overline{B(\overline
x_i,r_h)}\right] \label{holes} \\
U_h  & := \bigcup_{i\in I_h}\left[
B(x_i,R_h) \vee B(\overline x_i,R_h)\right].
	\label{annuli} 
\end{align}

\begin{definition}\label{h-relaxed-mfld}
For  $h \in \n$ we let 
$$
M_h : = S^d \setminus E_h,
$$ 
so that $\partial M_h = \partial E_h$, 
and let $g_h$  denote the metric on $M_h$ which is defined as the
standard metric $g$ of $S^d$ on $S^d\setminus U$, and the 
metric of the handle $(H_h,\varphi_{i,h}) $ on $ \Bigl( B(x_i,R_h)
\setminus  B(x_i ,r_h) \Bigr)  \vee   \Bigl( B(\overline x_i , R_h)
\setminus  B( \overline x_i ,r_h)  \Bigr) $, $i\in I_h$. Finally we
let $T$ be the antipodal map (\ref{antipodal-map}).
\end{definition}

\begin{remark}\label{h-curvature}
Notice that one of the sectional curvatures of $M_h$ is equal to 
$1/\e_h$. In particular the sectional curvature of $M_h$ it is
unbounded, as $h \uparrow +\infty$. 
\end{remark}

\begin{lemma}\label{limit-holes-annuli}
The sequence  $(1_{M_h} )_h $ converges in measure $\mbox{Vol}_g$  to
the constant  function  1; moreover   the sequence $(g_{ij,h}1_{M_h}
)_h$ converges in measure  $\mbox{Vol}_g$  to $ g_{ij}$, 
$i,j=1,\ldots, d$. 
\end{lemma}

\begin{pf}
Let $(x_i)_{i \in I_h}$ be the $\eta_h$-packing introduced above; let
moreover $(U_h)_h$ be as in (\ref{annuli}) above. To prove the lemma
it   is sufficient to show that $\mbox{Vol}_g(U_h)$ tends to zero, as
$h\uparrow +\infty$, and to this aim, we first prove an estimate on
$\# (I_h)$,  and then an estimate on the volume of geodesic balls.  
By Definition~\ref{packings}-(p$_1$), we get that
$$
\mbox{Vol}_g\left(\bigvee_{i \in I_h} B(x_i, \eta_h ) \right) = 
\sum_{ i \in I_h } \mbox{Vol}_g B( x_i, \eta_h ) \le
\mbox{Vol}_g(S^d).  
$$
Thus
\begin{equation}
\#( I_h ) \le \frac{ \mbox{Vol}_g(S^d) }{ \min_{ i \in I_h } \mbox{Vol}_g( B ( x_i, \eta_h )) }.
\label{estimate-packing}
\end{equation}
The measure  $\mbox{Vol}_g$ of geodesic balls $ B ( x_i, \eta_h )) $
can be computed by means of the following formula, which is a
particular case of a Bishop-type inequality:
\begin{equation} 
\frac{
\mbox{Vol}_g ( B ( x_i, \eta_h)  ) }{ \omega_d~ \eta_h^d  }= \left( \frac { \sin (\eta_h /\e_h) } {  (\eta_h /\e_h) }  \right)^{d-1} 
\label{bishop-inequality}
\end{equation}
for $i\in I_h$,  where $\omega_d$ is the $d$-dimensional euclidean
volume of the unit sphere.  Notice that,  as $\lim_{h \uparrow
+\infty} \eta_h / \e_h = 0$ (cf. \ref{pinco}),
the ratio at the left-hand side of   (\ref{bishop-inequality}) tends
to 1  as $h \uparrow +\infty$. 
Let us estimate the volume of $U_h $:
$$
\mbox{Vol}_g ( U_h  ) \le  \sum_{ i \in I_h } \mbox{Vol}_g ( B( x_i, R_h) ).
$$
By means of   (\ref{bishop-inequality}) and  (\ref{estimate-packing}), we have
$$ 
\mbox{Vol}_g ( U_h  ) \le \zeta_h  \e_h^d,
$$
where $(\zeta_h)_h$ is a bounded sequence. Thus, passing to the limit as $h \uparrow +\infty$, 
$ \mbox{Vol}_g ( U_h  )  $ tends to zero and  we get  the result.
\end{pf}

\begin{definition}\label{relaxed-funct}
For $h\in \n$, let ${\cal D}_{M_h}(\cdot)$ be the Dirichlet functional
on $ ( M_h, g_h)$, and  $f   \in L^2 (S^d ,  g)$; for  $\lambda > 0$,
let $ {\cal R}_h^\lambda : L^2 (S^d 
(1) , g ) \longrightarrow  [0, +\infty]$ be defined by 
\begin{align*}
{\cal R}_h ^\lambda (u) &  : = {\cal D}_{M_h} ( u ) + \int_{\overline
M_h} \bigl[ u ( x ) - u ( T ( x ) ) \bigr]^2 \infty_{\partial M_h} (
dx ) \\ 
& + \lambda \int_{M_h} u^2\mbox{Vol}_{g_h} (dx)  - 2\int_{M_h} f u~
\mbox{Vol}_{g_h} (dx) 
\end{align*}
if $ u|_{ M_h } \in H^1 (M_h, g_h)$; ${\cal R}_h (u) := +\infty$
otherwise in $L^2 (S^d(1), g)$. 
\end{definition}

\begin{prop}\label{h-representation}
(i)\quad The manifold \( (N_h,g^{N_h}) \) is
represented (in the sense of Definition~\ref{relaxriem}) by  \(
(\overline M_h,g_h,T,\infty_{\partial M_h}) \), $h\in \n$.  

(ii)\quad The sequence of manifolds $ (( M_h , g_h ))_h $ satisfies the
the following conditions:
\begin{itemize}
\item[(a$_1$)] There exists $\Lambda_o>0$ such that for every
$h\in\n$ 
\begin{equation}
0\le\Lambda_o^{-1}\lsumgxi\le\lsumghxi\le\Lambda_o\lsumgxi,
\label{comparemetric-1}
\end{equation}
for every $x\in M_h$ and every $\xi\in{\r}^d$;
\item[(a$_2$)] {\em (Strong Connectivity Condition)\/} There exist
bounded linear extension operators $\pi_h:\huh\longrightarrow\hu$
which satisfy the uniform bound
\[
\|\pi_hu\|_\hu \le c_o\|u\|_\huh,
\]
for every $u\in\huh$, and the constant $c_o$ does not depend on $h\in\n$.
\item[(a$_3$)] The sequence of characteristic
functions $(1_{M_h})_h$ converges in measure to the constant function
equal to 1 on $S^d$.
\item[(a$_4$)] \( T(M_h) = M_h \); 
\item[(a$_5$)] the canonical measure \( \mbox{Vol}_{g_h} (dx)=
\sqrt{\det{g_h }\,}~dx\) is \( T \)-invariant, \ie
\[
\sqrt{\det{g_h }\,}(x) = \sqrt{\det{g_h }\,}(T(x))\det\left(
\frac{\partial T}{\partial x}(x) \right),
\]
for almost every \( x\in M_h\); 
\item[(a$_6$)] the metric \( g_h \) is \( T \)-invariant, \ie
\[
g^{ij}_h(T(x)) = \sum_{k,\ell =1}^d \frac{\partial T^i}{\partial
x_\ell}(x) \frac{\partial T^j}{\partial x_k}(x) g^{\ell k}_h(x),
\]
for almost every \( x\in M_h\).
\end{itemize}
\end{prop}

\begin{pf} Arguing as in Proposition~\ref{representation} we can show
part (i).  Applying \cite[Proposition 3.9]{dmgm},
Definition~\ref{h-relaxed-mfld} and Lemma~\ref{limit-holes-annuli}, we
get (ii); more precisely,  from Proposition 3.9 in \cite{dmgm} we get
(a$_1$) and (a$_2$); (a$_3$) follows from
Lemma~\ref{limit-holes-annuli}; finally, by definition, $T ( M_h ) =
M_h$, both the  metric $g_h$ and the canonical measure
$\mbox{Vol}_{g_h}(\cdot)$ are invariant under the antipodal map $T$,
and moreover $\infty_{ \partial 
M_h}( \cdot ) = \infty_{ \partial M_h}( T ( \cdot ) )$, i.e.,
(a$_4$), (a$_5$), (a$_6$) are
satisfied. This proves the proposition. 
\end{pf} 

\begin{corollary}\label{rem-conv-spaces}
The sequence of Hilbert spaces $(L^2(N_h , g^{N_h}))_h$ is uniformly 
embedded into $L^2(S^d , g)$.
\end{corollary}

\begin{pf*}{Proof of Theorem~\ref{limit-handles}} 
Adapting  some arguments from  Proposition 5.7 in \cite{dmgm}, the proof
of the theorem follows if we show that the $\Gamma$-limit (in
$L^2( S^d, g)$) of the sequence $({\cal R}_h^\lambda)$ is equal to 
the functional
\begin{align*}
{\cal D}_\infty (u) & : = {\cal D}_{S^d} (u) + \frac{\alpha}{2^d}
\int_{S^d} \bigl[ u(x) - u(T(x)) \bigr]^2\mbox{Vol}_g(dx) \\
& \ \ + \lambda \int_{S^d} u^2 \mbox{Vol}_g(dx) - 2 \int_{S^d}
f u~ \mbox{Vol}_g(dx)
\end{align*}
for $u\in H^1(S^d , g)$. By a general result in
$\Gamma$-convergence (\cite{dm1}) this is equivalent to prove that the
sequence of functionals $G_h : L^2 (S^d , g) \longrightarrow
[0,+\infty] $ defined by 
$$
G_h(u) : {\cal D}_{M_h} (u) + \int_{\overline M_h} \bigl[ u(x) -
u(T(x)) \bigr]^2 \infty_{\partial M_h}(dx) 
$$
if $u|_{M_h} \in H^1(M_h g_h ) $, $G_h(u) = +\infty$ otherwise in $
L^2 (S^d , g)$, $\Gamma$-converges to the functional 
$$
G(u) : = {\cal D}_{S^d} (u) + \frac{\alpha}{2^d}
\int_{S^d} \bigl[ u(x) - u(T(x)) \bigr]^2\mbox{Vol}_g(dx),
$$
for $u\in H^1(S^d , g)$.

By Proposition~\ref{h-representation}-(ii) the sequence of manifolds
$(M_h , g_h) $ satisfies the assumptions (a$_1$), (a$_2$), (a$_3$),
(a$_4$), (a$_5$), (a$_6$), hence we can apply Theorems~\ref{var-comp} and
~\ref{handles-comp} in \S~\ref{varcompmflds} and get that  the
$\Gamma$-limit of $(G_h)$ is equal to 
$$
\widetilde G (u) = {\cal D}_{S^d} (u) +
\int_{S^d} \bigl[ u(x) - u(T(x)) \bigr]^2 \mu(dx),
$$
for some measure in ${\cal M}_o (S^d , g)$. What is left to prove,
then, is that $\mu(dx) = \alpha 2^{-d}\mbox{Vol}_g (dx)$; with our
choice of $\e_h$, plus the assymption (\ref{hypo-radii}) in
Theorem~\ref{limit-handles}, this can be done suitably modifying a
derivation-type argument as in \cite{bdmm2}. The proof of the theorem
is then complete. 
\end{pf*}

\section{Appendix}\label{varcompmflds} 
\noindent
For each $h\in\n$, let $(\overline M_h,g_h)$, be a
manifold-with-boundary, $\overline M_h =
M_h \cup \partial M_h$, $\partial M_h \not=\emptyset$, with
$\overline M_h\subset M$, where \(M \) is a manifold (with or without
boundary); we assume  that \( \mbox{dim}~M_h = d = \mbox{dim}~M\).

\subsection*{Variational compactness of Lipschitz metrics}
We shall consider the following assumptions in the rest of the paper. 
\begin{itemize}
\item[(A$_1$)] There exists $\Lambda_o>0$ such that for every
$h\in\n$ 
\begin{equation}
0\le\Lambda_o^{-1}\lsumgxi\le\lsumghxi\le\Lambda_o\lsumgxi,
\label{comparemetric}
\end{equation}
for every $x\in M_h$ and every $\xi\in{\r}^d$;
\item[(A$_2$)] {\em (Strong Connectivity Condition)\/} There exist
bounded linear extension operators $\pi_h:\huh\longrightarrow\hu$
which satisfy the uniform bound
\[
\|\pi_hu\|_\hu \le c_o\|u\|_\huh,
\]
for every $u\in\huh$. 
\newline
The constant $c_o$ does {\em not\/} depend on $h\in\n$.
\item[(A$_3$)] The sequence of characteristic
functions $(1_{M_h})_h$ converges in the weak$^*$ topology of
$L^\infty(M)$ to the function 
$b$; moreover {\em both \/}  \( b \) {\em and\/} \( b^{-1}\) belong to
\( L^\infty(M)\). 
\end{itemize}

\begin{remark}\label{volumecomparison}
1) \quad From (\ref{comparemetric}) in (A$_1$) we get the following
formula relating the local densities 
\begin{equation}
\Lambda_o^{-d/2}\locvol\le \locvolh \le \Lambda_o^{d/2}\locvol
\label{comparevolume}
\end{equation}
for every $x\in M_h$ and $h\in\n$. In particular sets of \(
\mbox{Vol}_{g_h}\)-measure zero are also of \(
\mbox{Vol}_{g}\)-measure zero, and conversely.

2) \quad Let $d(\cdot,\cdot)$, $d_h(\cdot,\cdot)$ be the distances on
$M_h$ associated with the metrics $ g $ (restricted to $M_h$),
$ g_h $ respectively. Then the assumption
(\ref{comparemetric}) implies that, for every $h\in\n$, the metric
space $(M_h,d_h)$ is equivalent to $(M_h,d)$.

3) \quad Note that from (\ref{comparemetric}) sets of \( (M_h,
g_h)\)-capacity zero has also \( (M,g) \)-capacity zero. Conversely,
if \( Z\subset M_h \) has \( (M,g)\)-capacity zero, then it also has
\( (M_h, g_h)\)-capacity zero.
\end{remark}

\vspace{.125in}
We are in a position to prove the following result.

\begin{th}\label{var-comp}
Let us assume (A$_1$), (A$_2$), (A$_3$), and consider
the sequence $({\cal D}_{M_h})_h$, where 
${\cal D}_{M_h}$ is defined by
\[
{\cal D}_{M_h}(u) := \left\{
\begin{array}{ll}
\displaystyle{ \int_M \sum_{i,j=1}^d g_h^{ij} D_i u D_j
u~\mbox{dVol}_{g_h}, \ \mbox{if}\ u|_{M_h} \in H^1(M_h,g_h),} \\
\\
\displaystyle{
+ \infty, \ \mbox{otherwise in}\ L^2(M,g),}
\end{array}\right.
\]
for \( h\in \n \). Then there exists a Lipschitz metric \(
a = \riema \) such that the sequence \( ({\cal D}_{M_h}) \)  
\( \Gamma\)-converges in \( L^2(M,g) \) to the weighted Dirichlet
functional  
\[
{\cal D}_{M,a,w}(u):=\left\{
\begin{array}{ll} 
\displaystyle{\int_{M_h}\lsumadu\sqrt{ \det a\,}~w(x)dx, \ \ 
u\in H^1(M,a),} & \\
\\
+\infty, \mbox{otherwise in}\ L^2(M,a), &
\end{array}
\right.
\]
with \( w(x) := \sqrt { \det g/\det a\, }(x) \) for  \( x\in M\), and
$\riema$ satisfies
\begin{equation}
\Lambda_o^{-1}\lsumgxi\le\lsumaxi\le c^2_o\lsumgxi,
\label{matrix-comp}
\end{equation}
for every $x\in M$ and every $\xi\in{\r}^d$; the constants $\Lambda_o$,
$c_o$ are those appearing respectively in (A$_1$), (A$_2$).
\end{th}

\begin{remark}\label{weight}
1) Notice that because of (\ref{matrix-comp}) both \(
w(\cdot) \) and \( w^{-1}(\cdot) \) are contained in \(
L^\infty(M)\). 

As another  consequence of (\ref{matrix-comp}) we have that if 
\( u\in H^1(M,g) \), then \( {\cal D}_{M,a,w}(u)<+\infty\); also, if 
\( v\in  H^1(M,a) \), then \( {\cal D}_M(v) < +\infty\). 

2) Using a similar argument as in Remark~\ref{volumecomparison}-3), we
have that a set has \( (M,a) \)-capacity zero if and only if has  
\( (M,g) \)-capacity zero. In particular, \( {\cal M}_o(M,a) \)
coincide with \( {\cal M}_o(M,g) \). 

3) From (\ref{matrix-comp}), and similarly as 
in Remark~\ref{volumecomparison}-1), we have that sets of \(
\mbox{Vol}_a\)-measure zero are also of \(
\mbox{Vol}_g\)-measure zero, and conversely.

4) Another consequence of (\ref{comparevolume}) and
(\ref{matrix-comp}) is that 
\[
\lim_{h\uparrow \infty} \|u_h - u\|_{L^2(M,g)} = 0 
\ \mbox{if and only if}\ 
\lim_{h\uparrow \infty} \|u_h - u\|_{L^2(M,a)} = 0 
\]
\end{remark}
\noindent
For the proof of Theorem~\ref{var-comp} we shall need the following
generalization of Theorem 2.1 by P. Marcellini \& C. Sbordone
\cite{MS}.
\begin{prop}\label{marsbo}
For each relatively compact open set $\Omega\subset M$ let us consider
for every 
$h\in\n$ the following functional 
\[
F_h(u,\Omega):=\left\{
\begin{array}{ll}
\displaystyle{\int_{M_h\cap\Omega}\usumghdu\locvolh~dx,
u\in \mbox{Lip}(M),} & \\ 
\\
+\infty, \mbox{otherwise in}\ L^2(M,g). &
\end{array}
\right.
\]
Then there exist a symmetric tensor $\uriema$ on $M$ and a functional
$F(\cdot,\Omega):L^2(\Omega, g)\longrightarrow [0,+\infty]$ such that 
\[
F(u,\Omega) = \int_\Omega \usumadu\locvol~dx,
\]
with $u\in \mbox{Lip}(M)$, and the tensor $\uriema$ satisfies
\[
0\le \usumaxi \le \Lambda_o\usumgxi,
\]
for all $x\in M$,  $\xi\in {\bf R}^d$ and $a^{ij} \in L^\infty(M)$, 
$i,j=1,\ldots,d$.
\end{prop}

\noindent
During the proof of the Proposition~\ref{marsbo} we shall need the
following general result \cite[Lemma 1.9]{n}.

\begin{lemma}\label{ln} 
Given any cover $(W_\ell)_{\ell\in I}$, $I\subseteq\n$,
of a paracompact, differentiable manifold $X$, and any Borel measure
$\sigma$ on $X$, there exists a family of open sets $(U_\ell)_{\ell\in
I}$ such that 
\begin{itemize}
\item[{\em (1)}] $\displaystyle{\bigcup_{\ell\in I} U_\ell = X \setminus
\left[~\bigcup_{\ell\in I}\partial U_\ell~\right]}$;
\item[{\em (2)}] $U_\ell \cap U_k = \emptyset$, for all $k,\ell\in I$, with
$k\not=\ell$;
\item[{\em (3)}] $U_\ell \subset \subset W_\ell$, with $\ell\in I$;
\item[{\em (4)}] $\displaystyle{\sigma\left[~\bigcup_{\ell\in I}\partial
U_\ell~\right]=0}$.
\end{itemize}
\end{lemma}

\begin{pf*}{Proof of Proposition~\ref{marsbo}} 
For each relatively compact open set \( \Omega \) in $M$, let
$F(\cdot,\Omega)$ be the functional which is the 
$\Gamma$-limit (in $L^2(M,g)$) of the sequence $(F_h)_h$. With
suitable modifications in  Lemmas 2.2 through 2.6 in \cite{MS} we can
prove that for each $u\in \mbox{Lip}(M)$ the set function $\Omega\mapsto
F(u,\Omega)$ satisfies the following properties: 
\begin{itemize}
\item for every relatively compact open sets
$\Omega'\subset\Omega\subset M$ 
\begin{align}
0 & \le  F(u,\Omega') ~ \le~  F(u,\Omega) 
\hphantom{~\le~F(u,\Omega')}  \nonumber \\ 
  \label{radon-nikodym} \\
& \le  F(u,\Omega') + \Lambda_o^{1+d/2} \int_{\Omega\setminus \Omega'}
\usumgdu\locvol~dx; \notag
\end{align}
\item for every disjoint relatively compact open sets $\Omega$,
$\Omega'$ in $M$ 
\begin{equation}
F(u,\Omega\cup\Omega') = F(u,\Omega) + F(u,\Omega').
\label{additivity}
\end{equation}
\end{itemize}
Using standard argument in Measure Theory we can extend the function 
$\Omega\mapsto F(u,\Omega)$ to a Borel measure 
\( \tau(u,\cdot) \) on $M$ with 
\[ 
\tau(u,\cdot) = F(u,\cdot)
\] 
on relatively compact open sets; moreover using 
a Radon-Nikodym argument, similarly as in the proof  \cite[Lemma 2.8]{MS}, we can show that $\tau(u,\cdot)$ is absolutely continuous w.r.t. $\mbox{Vol}_g(\cdot)$. If in the above
Lemma~\ref{ln} we let $(W_\ell)_{\ell\in I}$ be the given atlas of
$M$, and  $\sigma(\cdot)=\mbox{Vol}_g(\cdot)$ then we can apply Lemma
2.8 in \cite{MS} on each relatively compact open set
\[
\Omega\cap U_\ell, \ \ \ell\in I,
\] 
and find a matrix $(\alpha^{ij}_\ell)_{i,j=1}^d$ such that:
\[
0\le \sum_{i,j=1}^d \alpha_\ell^{ij}\xi_i\xi_j \le
\Lambda_0^{1+d/2}\locvol\sum_{i,j=1}^d g^{ij}\xi_i\xi_j,
\]
for all $x\in U_\ell\cap \Omega$, and $\xi\in\r^d$,
\[
\alpha_\ell^{ij} = \alpha_\ell^{ji}, \ \ \alpha_\ell^{ij}\in
L^\infty(M), i,j=1,\ldots,d,
\]
so that if we define
\[
a_\ell^{ij} := \frac{\alpha_\ell^{ij}}{\Lambda_o^{d/2}\locvol},
i,j=1,\ldots, d,
\]
we get
\begin{equation}
F(u,\Omega\cap U_\ell) = \int_{\Omega\cap U_\ell}
\sum_{i,j=1}^d a_\ell^{ij}D_iuD_ju\locvol~dx,
\label{rappell}
\end{equation}
for $u\in \mbox{Lip} (M)$. 
\newline
As the measure $\tau(u,\cdot)$ is absolutely continuous w.r.t.
$\mbox{Vol}_g(\cdot)=\sigma(\cdot)$,
we have from Lemma~\ref{ln}-(4) that
$\tau(u,\displaystyle{\bigcup_{\ell\in 
I}\partial U_\ell})=0$;   the  sets $U_\ell$, $\ell \in I$,  are disjoint, hence
\begin{math}
\displaystyle{
\tau(u,\Omega) = \sum_{\ell\in I}\tau(u,\Omega\cap U_\ell),}
\end{math}
which implies 
\begin{equation}
F(u,\Omega)=\sum_{\ell\in I}F(u,\Omega\cap U_\ell),
\label{measure}
\end{equation}
as  $\tau(u,\cdot)$ coincides with $F(u,\cdot)$ on relatively
compact  open
sets. Therefore by (\ref{rappell}) and (\ref{measure}) we get 
\begin{align*}
F(u,\Omega) & =  \sum_{\ell\in I}\int_{\Omega\cap U_\ell}
\sum_{i,j=1}^d a_\ell^{ij}D_iuD_ju\locvol~dx \\
& =  \int_\Omega \usumadu\locvol~dx,
\end{align*}
where the (0,2)-tensor $a=\uriema$ is   defined by 
\[
a^{ij}=a^{ij}(x) := \left\{
\begin{array}{ll}
a^{ij}_\ell(x), & x\in U_\ell \\
\\
g^{ij}(x), & x\in \bigcup_{\ell\in I}\partial U_\ell,
\end{array}
\right.
\]
so that we have
\[
0 \le \usumaxi \le \Lambda_o\usumgxi,
\]
for every $x\in M$ and for every $\xi\in{\bf R}^d$, with $a^{ij}=a^{ji}$,
$a^{ij}\in L^\infty(M)$, $i,j=1,\ldots,d$. The proof of the
proposition is thus completed.  
\end{pf*}

\begin{pf*}{Proof of Theorem~\ref{var-comp}}
Let $\phi:L^2(M,g)\longrightarrow [0,+\infty]$ be the $\Gamma$-limit
(in $L^2(M,g)$) of (a possible subsequence of)  $(\widetilde{{\cal D}}_{M_h, g_h}(u) )_h$. We notice that, for each $h\in\n$, the functional $(\widetilde{{\cal D}}_{M_h, g_h}(u) )_h$ is also  the lower semi-continuous regularization of the functional $F_h(\cdot,M)$
introduced in Proposition~\ref{marsbo}: this statement can be proven
as in  
the Step 1 of the proof of Theorem 4.4. in \cite{dmgm}. Hence it
is not difficult to see that  $ \phi$  coincides with $F(\cdot,M)$,
the $\Gamma$-limit of 
$(F_h(\cdot,M))_h$. The previous Proposition~\ref{marsbo}, then, gives
us the representation formula for $\phi$, namely, 
\[
\phi(u) = \int_M \usumadu\locvol~dx
\]
with $u\in \mbox{Lip} (M)$, and the tensor $\uriema$ satisfies 
\begin{equation}
0\le \usumaxi \le \Lambda_o\usumgxi,
\label{continuity}
\end{equation}
for all $x\in M$, and $\xi\in{\bf R}^d$, and
\begin{equation}
\left\{
\begin{array}{ll}
a^{ij}\in L^\infty(M), & \\
\\
a^{ij} = a^{ji}, &
\end{array}
\right.
\label{symmetry}
\end{equation}
for $i,j=1,\ldots,d$.
As in Step 2 of the proof of \cite[Theorem
4.4]{dmgm} we have
\[
\phi(u) = \int_M \usumadu\locvol~dx
\]
for all $u\in\hu$. 
\newline
We now prove that the functional $\phi$ is equal to $+\infty$ outside
$\hu$; more precisely, we have 
\[
\|u\|_\hu^2 \le c_o^2\phi(u) + c_0^2\Lambda_o^{d/2}
\|u\|_{L^2(M,g)}^2, 
\] 
and $\phi(u)=+\infty$ for every $u\in L^2(M,g)\setminus\hu$. 
\newline
First of all we notice  that 
\begin{equation}
\|v\|_{L^2(M_h,g_h)}^2 \le \Lambda_o^{d/2} \|v\|_{L^2(M,g)}^2.
\label{***}
\end{equation}
Indeed by (\ref{comparemetric}) we have
\begin{align*}
\|v\|_{L^2(M_h,g_h)}^2 & =  \int_{M_h}|v|^2\locvolh~dx
\le \Lambda_o^{d/2} \int_{M_h} |v|^2\locvol~dx \\
& \le  \int_M |v|^2\locvol~dx \le \Lambda_o^{d/2} 
\|v\|_{L^2(M,g)}^2.
\end{align*}
Let $u\in L^2(M,g)$ be such that $\phi(u)<+\infty$.  As $\phi$ is the
$\Gamma$-limit of $(\phi_h)_h$, there exists by  definition a sequence
$(u_h)_h$ of functions converging to  $u$ in $L^2(M,g)$ with
\[
\lim_{h\uparrow +\infty} \phi_h(u_h) = \phi(u).
\]
We notice that 
\[
\|u_h\|_\huh = \phi_h(u_h) + \|u_h\|_{L^2(M_h,g_h)}^2, 
\] 
and by (\ref{***})
\[
\|u_h\|_\huh \le \phi_h(u_h) + \Lambda_o^{d/2} \|u_h\|_{L^2(M,g)}^2. 
\]
Hence, by our assumption on $u$ and $(u_h)_h$, 
\[
\limsup_{h\uparrow +\infty} \|u_h\|_\huh 
\le \phi(u) + \Lambda_o^{d/2} \| u\|_{L^2(M,g)}^2 < +\infty.
\]
Using the strong connectivity condition (A$_2$) (we recall that $c_o$ is
independent of $h$)
\[
\|\pi_h u_h \|_\hu^2 \le c_o^2 \|u_h\|_\huh^2
\]
(and setting $v_h := \pi_hu_h$ for shortness) 
we have 
\[
\limsup_{h\uparrow +\infty} \|v_h\|_\hu^2 \le c_o^2\left(
\phi(u) + \Lambda_o^{d/2} \| u\|_{L^2(M,g)}^2\right).
\]
Therefore, up to a subsequence, $(v_h)_h$ converges to a function
$v\in \hu$ weakly in $\hu$ and, by Rellich's theorem, strongly in
$L^2(M,g)$. Using again (\ref{***}), applied this time to $u - u_h$ and $v - v_h$, we have
\begin{align}
\| u-u_h\|_{L^2(M_h,g_h)}^2 & \le  
\Lambda_o^{d/2} \|u-u_h\|_{L^2(M,g)}^2 \notag \\
\| v-v_h\|_{L^2(M_h,g_h)}^2 & \le 
\Lambda_o^{d/2} \|v-v_h\|_{L^2(M,g)}^2. \notag
\label{11}
\end{align}
Being $\pi_h$ an extension operator, we also have 
\[
v_h = u_h\ \ \mbox{(hence $u=v$) on}\ \ M_h.
\]
By the assumption (A$_3$), which we recall here,
\[
\left\{
\begin{array}{ll}
(1_{M_h})_h \ \mbox{converges $w^*$-$L^\infty(M)$ to}\ b, & \\
\\
b,b^{-1}\in L^\infty(M),
\end{array}
\right.
\]
and by H\"older's inequality we get
\begin{equation}
\| u-v\|_{L^2(M,g)}^2 \le \| b^{-1}\|_{L^\infty(M)}
\int_M |u-v|^2b(x)\locvol~dx.
\label{12}
\end{equation}
From  this inequality we get  that $u=v$ $\mbox{Vol}_g$-almost
everywhere on $M$; indeed by (\ref{12}), (A$_3$) and
(\ref{comparevolume}), we have
\begin{align*}
\int_M |u-v|^2\locvol~dx & \le  \|b^{-1}\|_{L^\infty(M)}
\int_M |u-v|^2b(x)\locvol~dx  \\
& \le  \|b^{-1}\|_{L^\infty(M)}\lim_{h\uparrow +\infty}
\int_{M_h} |u-v|^2\locvol~dx \\
& \le  \|b^{-1}\|_{L^\infty(M)} \Lambda_o^{d/2}
\lim_{h\uparrow +\infty} \int_{M_h} |u-v|^2\locvolh~dx \\
& =  0.
\end{align*}
Therefore $u=v$~~ $\mbox{Vol}_g$-almost everywhere on $M$, hence
\[
\phi(u)<+\infty \ \Longrightarrow u\in\hu,
\]
and 
\begin{align*}
\|u\|_\hu^2 = \|v\|_\hu^2 & \le
\liminf_{h \uparrow +\infty} \| v_h \|_{ H^1(M,g) } \\
& \le c_0^2\left(\phi(u) + \Lambda_o^{d/2} \|u\|_{L^2(M,g)}\right).
\end{align*}
The proof of the theorem will be achieved if we prove that the
tensor $\uriema$ satisfies the following ellipticity condition
\begin{equation}
c_o^{-2}\usumgxi \le \usumaxi.
\label{ellipticity}
\end{equation}
Indeed if (\ref{ellipticity}) holds true, then it follows that
$\uriema$ is invertible and its inverse $a=\riema$ is symmetric, being
$\uriema$ such, and $a_{ij}\in L^\infty(M)$,
$i,j=1,\ldots,d$. Thus $\bigl(M,a\bigr)$ is a Lipschitz
manifold, and from (\ref{continuity}) and (\ref{ellipticity}) we
get 
\[
\Lambda_o^{-1} \lsumgxi \le \lsumaxi \le c_o^2\lsumgxi
\] 
for every $x\in M$ and every $\xi\in{\r}^d$. 

Now we prove (\ref{ellipticity}) and to this aim we make use of the
family $(U_\ell)_{\ell\in I}$ of open sets as in Lemma~\ref{ln}
with $X=M$, $(W_\ell)_{\ell\in I}$ is the given atlas of $M$ and
$\sigma(\cdot)=\mbox{Vol}_g(\cdot)$. Consider, for $\ell\in I$, the
open set 
$U_\ell$, and let  $\psi\in \mbox{Lip}(M)$. In local coordinates, and making use of the same computation done in  Step 4 in
the proof of \cite[Theorem 4.4]{dmgm},  we then have
\begin{equation}
\begin{array}{ll}
\displaystyle{c_o^{-2}\int_{U_\ell}\Bigl(\usumgxi\Bigr)\psi^2\locvol~dx}
&   \\  
\hspace{.6in} \displaystyle{\le
\int_{U_\ell}\Bigl(\usumaxi\Bigr)\psi^2\locvol~dx}
\label{**}
\end{array}
\end{equation}
for every $\xi\in{\r}^d$.  We apply Lemma~\ref{ln} with $\sigma(\cdot) = \mbox{Vol}_g(\cdot)$ and get, for every $\psi\in \mbox{Lip} (M)$, 
\[ 
\begin{array}{ll}
\displaystyle{c_o^{-2}\int_M \Bigl(\usumgxi\Bigr)\psi^2\locvol~dx}
&   \\  
\hspace{.6in} \displaystyle{\le
\int_M \Bigl(\usumaxi\Bigr)\psi^2\locvol~dx,}
\end{array}
\]
for every $\xi\in{\r}^d$. Therefore, up to a redefinition on a set of 
$\mbox{Vol}_g$-measure zero, we get (\ref{ellipticity}). Define then 
\( {\cal D}_{M,a,w}(\cdot) := \phi(\cdot)\), and the proof of the
theorem is complete.
\end{pf*}

\begin{prop}\label{cormeasure}
In addition to the assumptions (A$_1$), (A$_2$), (A$_3$) let us
suppose that
\begin{eqnarray}
1_{M_h} \longrightarrow b & & \mbox{in measure on $M$}
\label{meas1} \\
g_{ij,h} 1_{M_h} \longrightarrow a_{ij} & & 
\mbox{in measure on $M$,}\ i,j=1,\ldots,d.  \label{meas2}
\end{eqnarray}
Then the same conclusion of Theorem~\ref{var-comp}
holds and the Lipschitz metric is equal to
$a=(a_{ij})_{i,j}$, where the coefficients $a_{ij}$ are given by
(\ref{meas2}) above. 
\end{prop}

\noindent
\begin{pf} It follows the lines of the proof of Proposition 4.5 in
\cite{dmgm}.
\end{pf}

\subsection*{Variational limits of handles}
In this part we shall be concerned with the asymptotic behavior of
sequences of functionals $\phi_h : L^2 ( M, g) \longrightarrow [0 ,
+\infty ] $ of the type
\begin{equation}
\phi_h(u) := {\cal D}_{M_h}(u) + \int_{\overline M_h}
\Bigl[u(x) - u(T(x))\Bigr]^2\mu_h(dx),
\label{pert-diri-funct}
\end{equation}
if   \( u|_{M_h} \in H^1(M_h, g_h) \);  \(
\phi_h (u) := +\infty \) otherwise in \( L^2(M,g) \).  

In (\ref{pert-diri-funct}) above, \( (\mu_h)_h \) is a sequence of
measures such that \(  \mu_h \in {\cal M}_o(M_h,g_h) \), $h \in N$;
the map \( T: M\longrightarrow M \) is an isometry, \( T\circ T =
\mbox{id}_M\), and the fixed-point set of \( T\), \( \mbox{Fix}(T)\),
is a submanifold of \( M\).  

\begin{th}\label{handles-comp}
Assume that  \( \big( (M_h,g_h)\big)_h \) satisfies (A$_1$), (A$_2$),
(A$_3$) above and is a $T$-invariant sequence of Lipschitz manifolds,
\ie 
\begin{itemize}
\item[(A$_4$)] \( T(M_h) = M_h \); 
\item[(A$_5$)] \( \displaystyle{
\sqrt{\det{g_h }\,}(x) = \sqrt{\det{g_h }\,}(T(x))\det\left(
\frac{\partial T}{\partial x}(x) \right),
} \) a.e. in  \( M_h\); 
\item[(A$_6$)] 
\( \displaystyle{ 
g^{ij}_h(T(x)) = \sum_{k,\ell =1}^d \frac{\partial T^i}{\partial
x_\ell}(x) \frac{\partial T^j}{\partial x_k}(x) g^{\ell k}_h(x),
}\) a.e. in \(  M_h\).
\end{itemize}
Let us consider, for each \( h\in \n\), the relaxed  Dirichlet
functional $\phi_h(u)$ defined on $L^2(M,g)$ by 
\[
\phi_h(u) := {\cal D}_{M_h}(u) + \int_{\overline M_h}
\Bigl[u(x) - u(T(x))\Bigr]^2\mu_h(dx),
\]
if \( u\in L^2(M,g) \) with \( u|_{M_h} \in H^1(M_h, g_h) \); \(
\phi_h (u) := +\infty \) otherwise in \( L^2(M,g) \).
Assume that 
\( \mu_h \in {\cal M}_o(M_h) \), and \( \mu_h(\cdot) \sim
\mu_h(T^{-1}(\cdot)) \) (in the sense of Definition~\ref{mzero}).

Then there exists a \( T\)-invariant Lipschitz metric \( a=
(a_{ij})_{i,j}\) on \( M \), with a \( T\)-invariant canonical measure
\( \mbox{Vol}_a \), and a Borel measure \( \mu \in {\cal M}_o(M) \),
with \( \mu(\cdot) \sim \mu(T^{-1}(\cdot)) \), such that the sequence
of functionals given by (\ref{pert-diri-funct}) \( \Gamma\)-converges
in \( L^2(M,g) \) to the functional \( \phi(u)\) defined by
\begin{equation}
\phi(u) := {\cal D}_{M,a,w}(u) + \int_{\overline M}
\Bigl[u(x) - u(T(x))\Bigr]^2\mu(dx),
\label{pert-limit-funct}
\end{equation}
if \( u\in H^1(M,g) \), and \( \phi(u) := +\infty\) otherwise in \(
L^2(M,g) \). The functional \( {\cal D}_{M,a,w}(\cdot) \) is 
the weighted Dirichlet functional introduced in
Theorem~\ref{var-comp}. 
\end{th}

\begin{remark}\label{proofofth}
Using (A$_1$), (A$_2$), (A$_3$), the existence of a Lipschitz
metric follows from Theorem~\ref{var-comp},
where it is also proved that the sequence \( ({\cal D}_{M_h}(\cdot))
\) 
\(\Gamma\)-converges to \( {\cal D}_{M,a,w}(\cdot) \). 

The \( T\)-invariance of the canonical measure \( \mbox{Vol}_a\) and
of the metric \( a=(a_{ij})_{i,j} \) are proved similarly as in
\cite[Lemma 5.1]{dmgm}, by means of (A$_4$), (A$_5$), (A$_6$).

What is left to prove, then, is the existence of a Borel measure \(
\mu\in {\cal M}_o(M) \), with \( \mu(\cdot) \sim \mu(T^{-1}(\cdot))
\), such that \( (\phi_h) \) \(\Gamma\)-converges in \( L^2(M,g) \) to
(\ref{pert-limit-funct}). 
\end{remark}

\noindent
To this aim we need the following lemma, which is a generalization to
this framework of a result in the euclidean setting by
G. Buttazzo, G. Dal Maso, \& U. Mosco in \cite[Section 4 \&
Appendix]{bdmm1}; cf. also \cite[Lemma 5.3]{dmgm}.

\begin{lemma}\label{meas-recon}
Let \( m : H^1(M,a) \longrightarrow [0,+\infty] \) be such that:
\begin{itemize}
\item[(1)] If \( 0 \le u \le v \) a.e. on \( M\), then \( m(u) \le
m(v) \); 
\item[(2)] \( m(|u|) = m(u)\);
\item[(3)] \( m(u+v) \le m(u) + m(v) \), if \( \min\{ u(x) , v(x) \} =
0\) for almost every \( x\in M\); 
\item[(4)] \( m(u) = \lim_{n \to +\infty} m(u_n) \), for every
increasing sequence of positive functions \( (u_h)\) such that 
\[ 
(M,a)\mbox{-cap} \{ x\in M : u_h(x)\ \mbox{does not converge to}\
u(x) \} = 0;
\]
\item[(5)] \( m(0) = 0\); \( m(tu) = t^2m(u)\), for \(t\in\r\); \(
m(u+v) + m(u-v) = 2\bigl[m(u) + m(v)\bigr]\). 
\end{itemize}
Then there exists a measure \( \mu\in{\cal M}_o (M, a) \) such that 
\[
m(u) = \int_{\overline M} u^2d\mu,
\]
for every \( u\in H^1(M,a)\).
\end{lemma}

\begin{pf}It is a simple modification of Theorem 2.3
in \cite{n}, by means of \cite[Theorem 2.22]{n}.
\end{pf}

\noindent
{\em Proof of Theorem~\ref{handles-comp}.} By means of
Lemma~\ref{meas-recon}, the 
proof now runs parallel to the proof of Proposition 5.2 in
\cite{dmgm}, which gives the existence of a measure \( \mu\in {\cal
M}_o(M) \), with \( \mu(\cdot) \sim \mu(T^{-1}(\cdot))\), such that \(
(\phi_h) \) \( \Gamma\)-converges in \( L^2(M,g) \) to the functional
given in (\ref{pert-limit-funct}). The theorem is so proved. \qed

\bibliographystyle{plain}

\bibliography{references}

\end{document}